\documentclass[12Pt]{article}
\usepackage{amssymb,euscript,latexsym,amsmath}
\usepackage{graphicx}
\usepackage[dvips]{epsfig}
\usepackage{epstopdf}
\usepackage{setspace}
\usepackage{subfigure} 
\usepackage{color}

\usepackage{caption2}

\DeclareGraphicsExtensions{.eps,.ps,.eps}

\begin{document}

\title{\Large \bf  Multiply-interacting Vortex Streets}

\author{ \normalsize 
Babak G. Oskouei,  Eva Kanso and Paul K. Newton
\\[2ex]
{\footnotesize Aerospace and Mechanical Engineering, University of Southern California} \\
{\footnotesize 854 Downey Way, Los Angeles, CA 90089-1191} \\
 }

\date{\normalsize February 4, 2010}

\maketitle

\maketitle

{\small
\singlespacing
\begin{abstract}
We investigate the behavior of an infinite array of (reverse) von K\'{a}rm\'{a}n streets. 
Our primary motivation is to model the wake dynamics in large fish schools. 
We ignore the fish and focus on the dynamic interaction of multiple wakes where each wake 
is modeled as \textcolor{black}{a reverse} von K\'{a}rm\'{a}n street. 
There exist configurations where the infinite array of vortex streets is
in relative equilibrium, that is, the streets move together with the same translational velocity.  
We examine the topology of the streamline patterns in a frame moving
with the same translational velocity as the streets which lends insight into fluid 
transport through the mid-wake region.  Fluid is advected along different 
paths depending on the distance separating two adjacent streets. 
Generally, when the distance between the streets is large enough, 
each street behaves as a single von K\'{a}rm\'{a}n street and fluid \textcolor{black}{moves}
globally between two adjacent streets. When the streets get closer to each other, 
the number of streets that enter into partnership in transporting fluid among 
themselves increases. This observation motivates a bifurcation analysis 
which links the distance between streets to the maximum number of streets 
transporting fluid among themselves. 
We also show that for short times, the analysis of streamline 
topologies for the infinite arrays of streets can be expected
to set the pattern for the more realistic case of a finite array of truncated streets, 
\textcolor{black}{which is not in an equilibrium state and its dynamic evolution}
eventually destroys the exact topological patterns identified 
in the infinite array case. 
The problem of fluid transport between adjacent 
streets \textcolor{black}{may be relevant} for understanding the transport of oxygen and nutrients to inner 
fish in large schools as well as understanding flow barriers to passive locomotion.
\end{abstract}
}

%%%%%%%%%%%%%%%%%%%%%%%%%%%%%%%%%%%%%%%

\doublespacing

\section{Introduction}
\label{sec:intro}

This paper considers the interaction of multiple reverse von K\'{a}rm\'{a}n streets. The primary motivation is to model the wake dynamics in large fish schools and to gain insight into the role of the fluid in transporting oxygen and nutrients to inner fish as well as its role in facilitating or acting as flow barriers to passive locomotion.

\textcolor{black}{We consider the parallel translational motion observed in fish schools
and study the dynamic interaction of the multiple wakes. Each wake is modeled 
as a reverse von K\'{a}rm\'{a}n street
and the fish producing the wakes are ignored.}
This model is reminiscent of the one employed in~\cite{Weihs1973b,Weihs1975}
in analyzing the hydrodynamic advantages for fish schooling. The focus
of~\cite{Weihs1973b,Weihs1975} was on investigating the energy-optimal positioning 
of an individual fish within the school (the famous diamond pattern) whereas
in the present study, we focus on the interaction and fluid
transport among multiple vortex streets for its relevance to understanding transport
of nutrients as well as barriers to passive locomotion.

The classical von K\'{a}rm\'{a}n vortex street -- made up of two rows of 
evenly spaced point vortices of equal and opposite sign staggered with respect to each other, 
see~\cite{Karman1911}  --  is the canonical `idealized' vorticity configuration appearing 
generically in the wake of bluff bodies; see, for example,~\cite{Williamson1996}.
While the dynamics and stability of the single street is well-understood (see, for example,~\cite{ArBrSt2007,Saffman1992} and references therein), little work is done on the interaction of multiple streets. 
\textcolor{black}{An infinite array of streets can be viewed, under certain conditions
on street alignments, as a special case of doubly-periodic vortex lattices
whose Hamiltonian dynamics is investigated in~\cite{Oneil1989, StAr1999}; see also~\cite{Oneil2007, Stremler2003} for investigations of singly-periodic lattices.
Given the results in~\cite{StAr1999}, one can argue that there exist configurations 
where an infinite array of streets is in relative equilibrium.}
This means that the streamline pattern remains steady 
in the frame moving with the same translational velocity as the streets. 
We look at the topology of the streamline patterns, that is to say, the streamline patterns relative 
to the moving frame, which lends insight into fluid transport through the mid-wake region.  Fluid is 
advected along different paths depending on the distance separating two adjacent streets. Generally, 
when the distance between the streets is large enough, each street behaves as a single 
von K\'{a}rm\'{a}n street and fluid is transported globally between two adjacent streets. 
When the streets get closer to each other, the number of streets that enter into partnership in transporting fluid among themselves increases. This observation motivates a bifurcation analysis which links the distance between streets to the maximum 
number of streets transporting fluid among themselves.

An investigation of the stability of the relative equilibria in the infinite array of vortex streets  
is central to assessing the prevalence of the reported streamline topologies under spatial and 
temporal perturbations. The stability of a square lattice (of identical vortices) is addressed in 
the classical work of~\cite{Tkachenko1966a,Tkachenko1966b}.
\textcolor{black}{Instead of casting the infinite arrays of streets as a Hamiltonian system following the methods developed in~\cite{Oneil1989, Oneil2007, Stremler2003,StAr1999} and studying its stability subject to periodic perturbations,
we choose to examine the dynamics of a finite array of vortex streets. The finite array consists of vortex streets 
that} are truncated horizontally
(reflecting the fact that in real flows viscosity diminishes the presence of the entire wake) and vertically
(reflecting a finite number of fish within the school). This approach does not constitute
a rigorous stability analysis of the relative equilibria of the infinite lattice. In fact, the
truncated array is no longer in a state of equilibrium but it could be interpreted as an arbitrary
perturbation to the infinite array whose dynamic evolution leads to valuable insights
into the time scales for which the streamline topologies observed in the infinite case remain valid.
Indeed, our study shows that many of the features observed in the idealized infinite case
are present up to a finite time in the truncated case. It also nicely demonstrates the mechanism by which 
these features get destroyed.

The organization of this paper is as follows.
The problem of an infinite array of vortex streets 
is described in Section~\ref{sec:formulate}
where this infinite configuration is shown to be
in relative equilibrium. Streamline topologies
which support communication and fluid transport between
multiple adjacent streets are reported in Section~\ref{sec:stream}.
In Section~\ref{sec:bif}, we present
a full bifurcation analysis that links the distance between the streets
to the maximum number of streets in communication. A discussion of the 
streamline topologies and their time evolution in the more realistic case of a finite 
array of truncated streets is presented in Section~\ref{sec:timescales}. 
The findings of this work are summarized in Section~\ref{sec:conc}.

%%%%%%%%%%%%%%%%%%%%%%%%%%%%%%%%%%%%%%%%%%%%%%%%%%%%%%%%%
\section{Infinite array of vortex streets}
\label{sec:formulate}

By way of background, consider the classical von K\'{a}rm\'{a}n 
street made up of two rows of evenly spaced 
point vortices of equal and opposite strength $\pm\Gamma$.
Let $a$ denote the distance between two neighboring vortices 
of the same row and $b$ denote the distance between the two rows. 
We are particularly interested in the staggered configuration
which is the canonical `idealized' vorticity configuration appearing 
generically in the wake of bluff bodies.
Without loss of generality, let the vortices of strength $-\Gamma$ 
be placed at $z = ma$ and those of strength $\Gamma$ 
be placed at $z = (m+1/2)a + i b$, 
$(-\infty<m<\infty)$. The complex notation ${z} = {x} + i {y}$ 
(where $i= \sqrt{-1}$) is employed 
for convenience.
%-------------------
\begin{figure}
   \centering
\subfigure[fixed frame]
  {
   \includegraphics[width=6cm,angle=0] {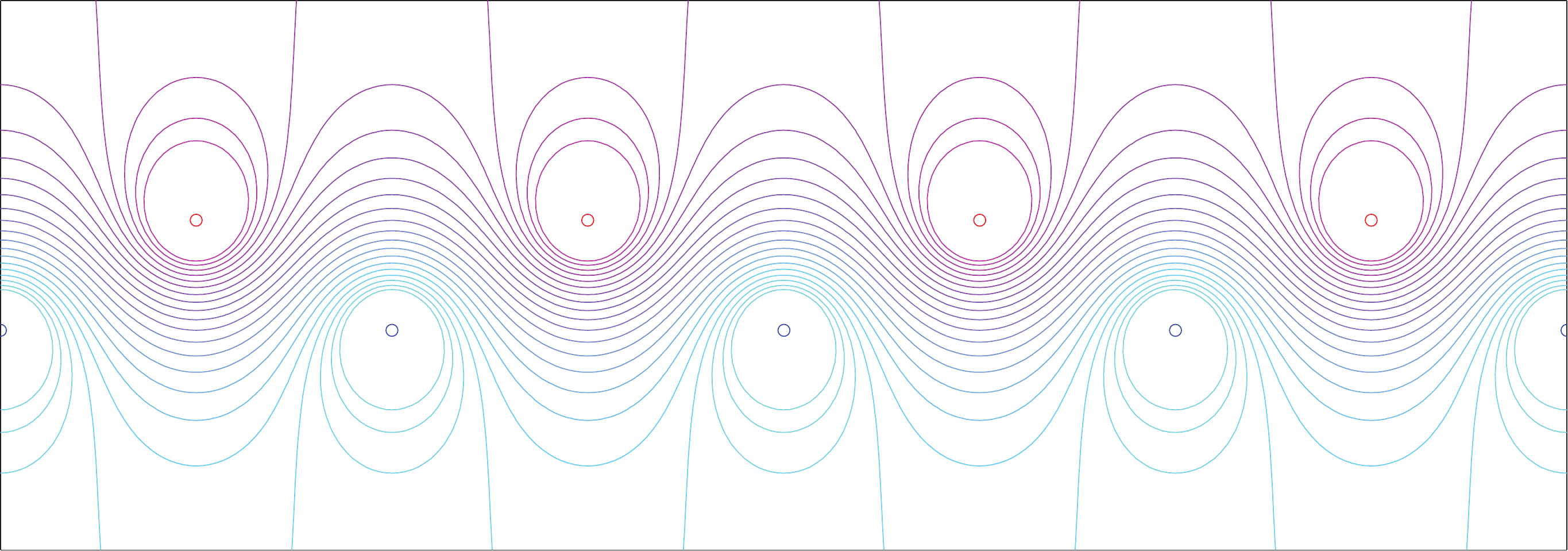}
  }
 \hspace{.5cm}
\subfigure[moving frame]
   {
   \includegraphics[width=6cm,angle=0] {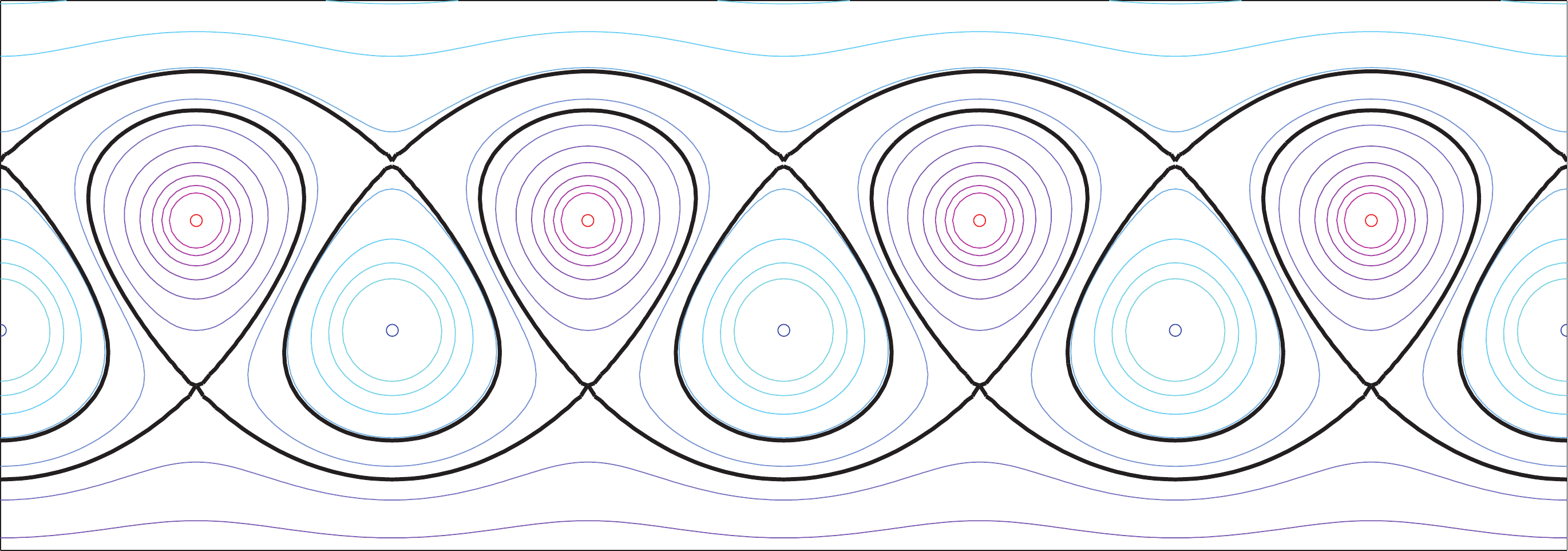}
   }
\caption{\textit{Staggered vortex street: (a)  streamlines plotted in fixed frame 
and (b) streamlines plotted in a frame moving the same translational velocity as the street, 
$U_{street} =  \frac{\Gamma}{2a} \tanh \frac{\pi b}{a}$. 
The parameter values are set to $\Gamma =1$, $a=1$ and  $b=0.2805$ which
corresponds to the value at which the staggered street is linearly stable. 
}}
\label{fig:singlestreet}
\end{figure}
%-------------------
The complex potential is found by superposition of the potential functions of all vortices 
$\lim_{M\rightarrow \infty} \sum_{m = - M}^{m = M}   \dfrac{\Gamma}{2\pi i}  
\log \dfrac{z - (m+1/2)a - i b}{z - ma}$
which converges to
%--------
\begin{equation}\label{eq:complex}
w_{street} = 
%\lim_{N\rightarrow \infty} \sum_{m = - N}^{m = N}   \dfrac{\Gamma}{2\pi i}  \log \dfrac{z - (m+1/2)a - i b}{z - ma}= 
\dfrac{\Gamma}{2\pi i} \log\left( \dfrac{\sin(\dfrac{\pi}{a}(z - a/2 - ib))}{\sin(\dfrac{\pi}{a} z)} \right),
\end{equation}
%--------
with $w_{street}\rightarrow 0$ as $y \rightarrow \pm \infty$; see, for example,~\cite{Saffman1992}
as well as the comprehensive discussion in~\cite{KoKiRo1964}(chapter 5, section 3). 
One can readily verify that this configuration corresponds to a relative equilibrium where
all the vortices move with the same translational velocity in the $x$-direction, namely,
%--------
\begin{equation}
{U}_{street} = \dfrac{\Gamma}{2a} \tanh \dfrac{\pi b}{a},  
\end{equation}
%--------
The staggered configuration can be shown to be neutrally stable (but nonlinearly unstable) for 
$\color{black}b/a = \frac{\cosh^{-1}\sqrt{2}}{\pi}\approx 0.2805$ and unstable
otherwise.
%\textcolor{black}{Despite the fact that the point vortex street model is unstable, physical effects such as viscosity allow vortex streets to exist in nature, which justifies the use of idealized point vortex street as model of a real wake.} 
In Figure~\ref{fig:singlestreet} is a depiction of the streamlines 
for the linearly stable case with $a=1$ and $\Gamma  = 1$. 
Figure~\ref{fig:singlestreet}(a) 
depicts the streamlines in inertial frame while Figure~\ref{fig:singlestreet}(b) 
shows the streamlines in a frame moving with the same translational velocity as the street itself.
More specifically, the streamlines in Figure~\ref{fig:singlestreet}(a) 
are given by level curves of the stream function  Im$(w_{street})$, the imaginary part of $w_{street}$, 
while those in Figure~\ref{fig:singlestreet}(b) 
correspond to level curves of Im$(w_{street}-  z U_{street})$, the imaginary part of $w_{street} -  z U_{street}$. 
The `thick' streamlines shown in Figure~\ref{fig:singlestreet}(b)
are particularly interesting. These streamlines are associated with 
stagnation points (in the moving frame) of hyperbolic or saddle type.
They are referred to as separatices because they 
separate regions with different fluid behavior. For example, one distinguishes
three regions in Figure~\ref{fig:singlestreet}(b).
In a region containing a point 
vortex and bounded by a separatix, fluid orbits around the vortex, whereas in the 
middle region bounded by both separatrices and void of point vortices, \textcolor{black}{the motion
of the fluid is global} while touring around a counterclockwise then 
a clockwise vortex. \textcolor{black}{Global motion} of the fluid is observed in 
the region bounded by one of the separatrices and the bound at infinity.

\begin{figure}
   \centering
   \includegraphics[scale = 0.8,angle=0] {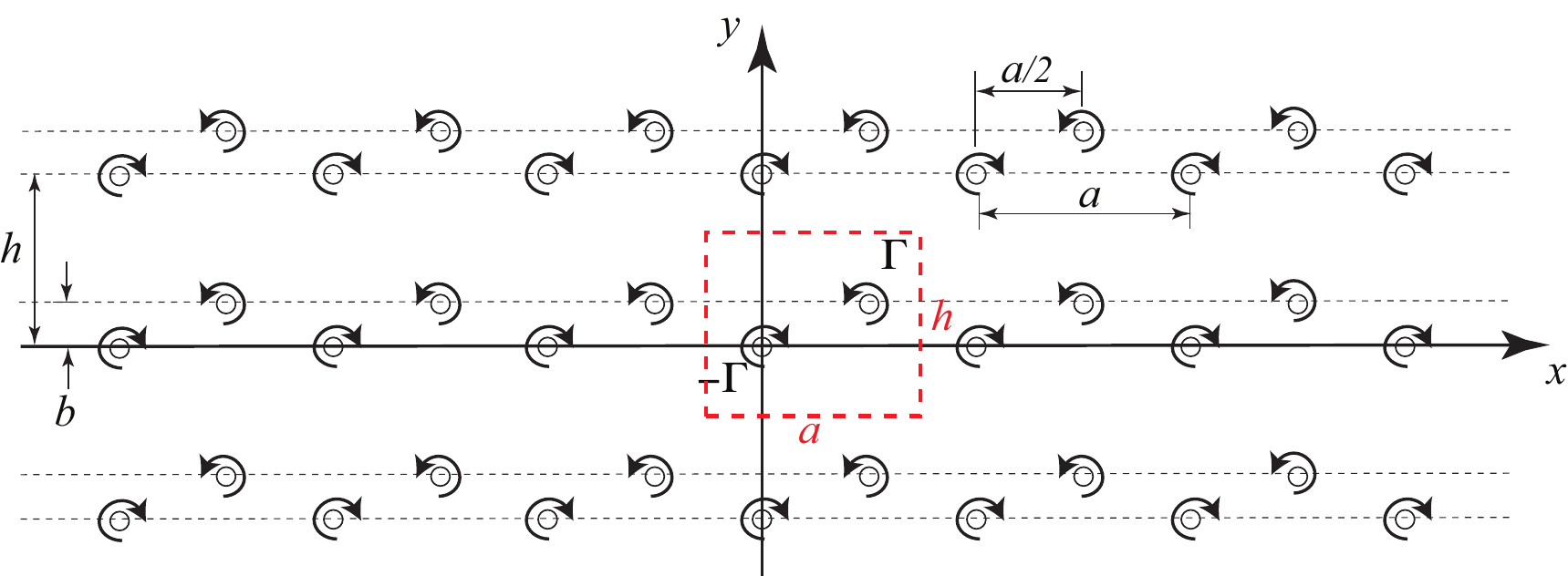}
\caption{Infinite array of staggered vortex streets: this configuration
can be viewed as a doubly-periodic vortex lattice or, as referred to in~\cite{StAr1999}, a 
vortex pair in a doubly-periodic parallelogram of length $a$ and width $h$.}
\label{fig:infinitestreets}
\end{figure}

We are interested in studying the configuration composed of an infinite array of vortex streets 
where each street has the same staggered structure, see Figure~\ref{fig:infinitestreets}. 
Let the vortices of strength $-\Gamma$ be located at $ma + nhi$  while the vortices of strength $\Gamma$ 
at $(m + 1/2)a +  (b + nh)i$, where $m$ and $n$ take the values $...,-2,-1,0,1,2,...$, and $n$ can be thought of
as identifying a whole street. 
%-------------------
%-------------------
\textcolor{black}{This configuration can be viewed as a vortex pair in a doubly-periodic parallelogram
(see Figure~\ref{fig:infinitestreets}) which is shown in~\cite{StAr1999} to be in relative equilibrium 
with the velocity of the vortices described in terms of the Weierstrass zeta-function. The fact
that this configuration is in relative equilibrium can be verified using a heuristic argument
which consists  of computing the velocities 
of two arbitrary vortices of opposite strength and showing that they are equal to a finite velocity $U$. 
From the symmetry of the infinite configuration, this is equivalent to showing that the 
velocity of the vortex of strength $-\Gamma$
placed at the origin $z = 0$ is equal to the velocity of, say, the vortex of strength $\Gamma$
placed at $z = a/2 + ib$.}

\textcolor{black}{We begin by writing the complex potential of the infinite array 
of streets, which is obtained  by superposition similarly to the single street case, namely,}
%-------------------
\begin{equation}\label{eq:doublesum}
w = \lim_{N\to\infty, M\to\infty}\sum_{n=-N}^{n=N}  \sum_{m=-M}^{m = M} 
 \dfrac{\Gamma}{2\pi i}  \log \dfrac{z - (m+1/2)a - i (b+nh)}{z - ma-inh} ,
\end{equation}
%-------------------
where $N$ and $M$  are integers ($N,$ $M \in \mathbb{N}$). Clearly, 
the inner sum corresponds to a sum over the $n^{\rm th}$ vortex street, which converges to
an expression similar to that in~\eqref{eq:complex}. 
Consequently, the complex potential~\eqref{eq:doublesum} can be written as
%-------------------
\begin{equation}\label{eq:infinitepotential}
w = \lim_{N\to \infty} w_N, \qquad w_N = \sum_{n=-N}^{N} \frac{\Gamma}{2{\pi}i} \log{\left[\frac{\sin{\frac{\pi}{a}\left(z-a/2-(b+nh)i\right)}}{\sin{\frac{\pi}{a}\left(z-nhi\right)}}\right]}.
\end{equation}
%-------------------
The complex velocity at a point that does not coincide with a point vortex 
is obtained by differentiating~\eqref{eq:infinitepotential} with respect to $z$. 
This yields
%-------------------
\begin{equation}\label{eq:complexvelocity}
\bar{u}=\lim_{N\to\infty}\sum_{n=-N}^{N}\frac{\Gamma}{2ai}\left[\cot{\frac{\pi}{a}
\left(z-a/2-(b+nh)i\right)}-\cot{\frac{\pi}{a}\left(z-nhi\right)}\right],
\end{equation}
%-------------------
where the bar is used to denote the complex conjugate (e.g., $\bar{z} = x - i y$). 
%-------------------
\begin{figure}[!t]
   \centering
   \includegraphics[scale = 0.75] {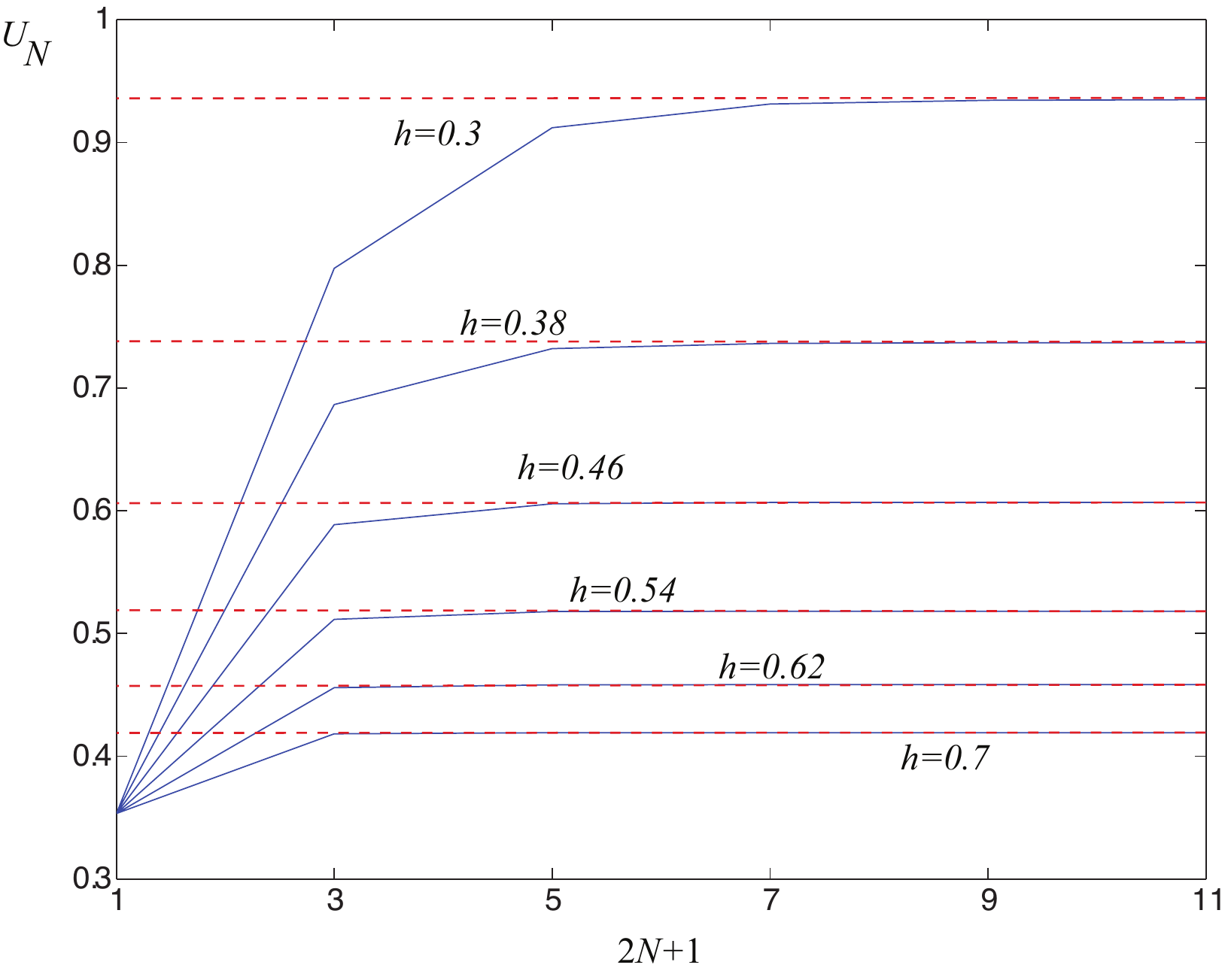}   
   \caption{\textit{Plot of the velocity $\color{black} U_N = \frac{\Gamma}{2a}\sum_{n=-N}^{N} \tanh(\pi(b + nh)/a)$ versus the total number of streets $2N+1$. 
   The parameter values are set to $a=1$, $\Gamma=1$, $b=0.2805$ and $h= 0.3, 0.38, 0.46, 0.54, 0.62, 0.7.$
   Similar convergence behavior is observed for various
values of $b<a/2$. The dashed lines correspond to the value of the equilibrium velocity obtained from the Weierstrass zeta-function based on the results in~\cite{StAr1999} plus the constant term $\Gamma b/ah$.}}
\label{fig:velocity}
\end{figure}
%-------------------
%\todo{Babak, could you please modify Figure 3 to include the he equilibrium velocity obtained from the Weierstrass zeta-function based on the results in~\cite{StAr1999} plus the constant term $\Gamma b/ah$}
The velocity of a point vortex is obtained by first subtracting the contribution
of that vortex from~\eqref{eq:doublesum} or~\eqref{eq:infinitepotential}, 
then differentiating the result with respect to $z$.
Actually, one could subtract the contribution of the whole row containing that vortex
because the effect of a given row on itself is always zero.\footnote{It is a known result that, for
an infinite row of equally spaced vortices of equal strength, all the vortices are at rest.}
To this end, the velocity of the vortex placed at $z = 0$ is given by
%-------------------
\begin{equation}\label{eq:z0velocity}
\bar{u}\vert_{z=0}=\lim_{N\to\infty}\dfrac{\Gamma}{2ai}
\left[ 
\sum_{n=-N}^{N} \cot  \dfrac{\pi(-a/2-(b+nh)i)}{a} \ \ -
\sum_{n=-N, n\neq 0}^{N}\cot \dfrac{-\pi nh i }{a}
\right] .
\end{equation}
%-------------------
One can readily verify that the second sum in~\eqref{eq:z0velocity} is zero\footnote{This follows from
writing $\sum_{n=-N, n\neq 0}^{N}\cot \dfrac{-\pi nh i }{a}= \sum_{n =1}^N \left(\cot \dfrac{-\pi nh i }{a} 
+ \cot \dfrac{\pi nh i }{a} \right) = 0$.}
while the first sum can be rewritten, using known trigonometric identities, 
in the  form
%------------------
\begin{equation}\label{eq:z0vel}
\bar{u}\vert_{z=0}=\lim_{N\to\infty}\dfrac{\Gamma}{2a}
\sum_{n=-N}^{N} \tanh \dfrac{\pi (b+nh)}{a} .
\end{equation}
%----------------
%---------------------
\begin{figure}
  \centering
   \subfigure[fixed frame]
   {
   \includegraphics[width=6cm,angle=0] {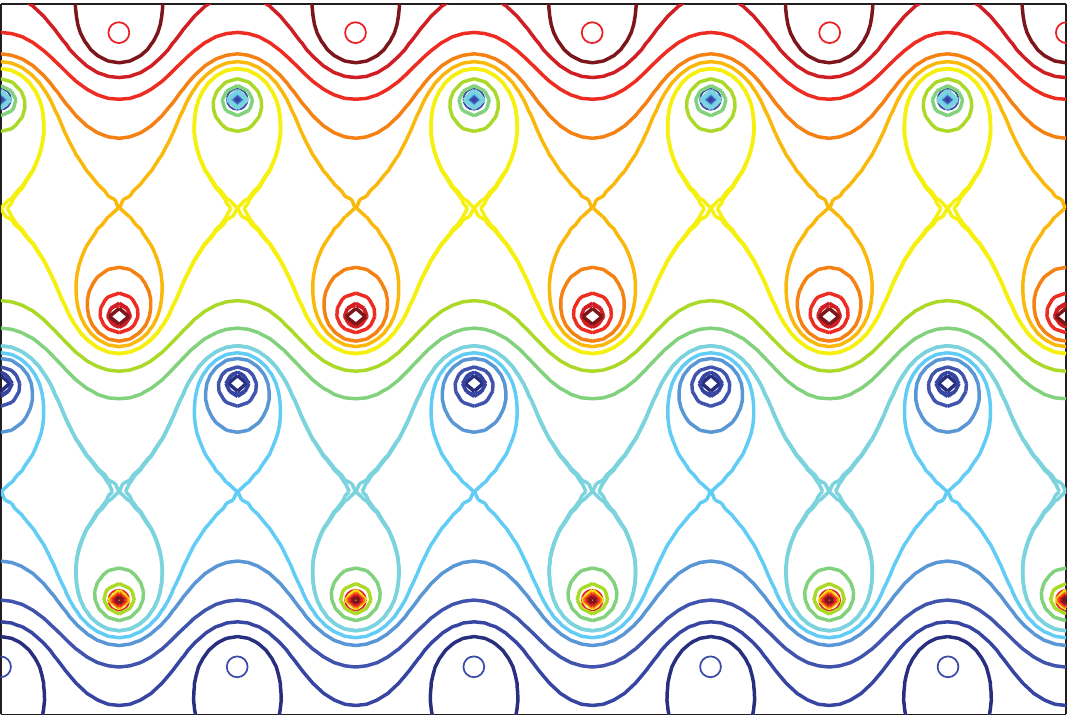}
   }
   \hspace{.5cm}
   \subfigure[moving frame]
   {
    \includegraphics[width=6cm,angle=0] {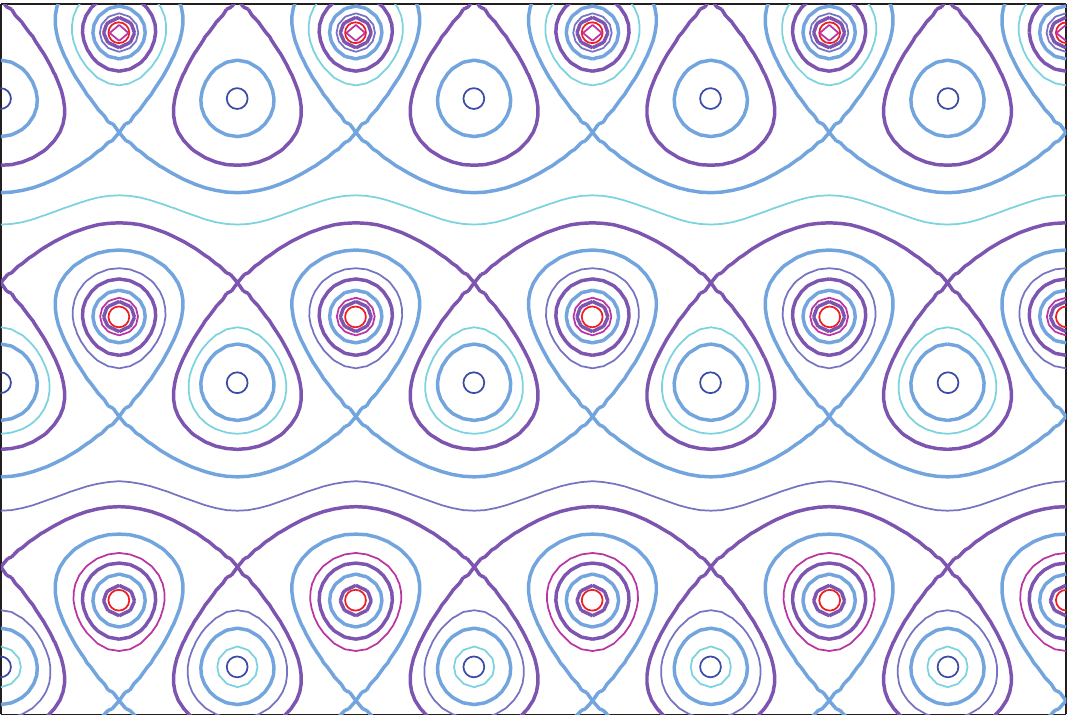}      
   }
   %%  \vspace{-1.5cm}
   \caption{\textit{Array of vortex streets: (a) streamlines plotted in fixed frame and (b) streamlines plotted in a frame
   moving with the same translational velocity as the streets. The parameter values are set to $\Gamma=1$, $a=1$, $\color{black} b=\frac{\cosh^{-1}\sqrt{2}}{\pi}$ and $h=1.2$. The infinite sums in~\eqref{eq:infinitepotential} and~\eqref{eq:velocity} 
 are truncated using $N=150$, i.e., using a total of $301$ streets.}}
\label{fig:far}
\end{figure}
%---------------------
Similarly, the velocity of the vortex placed at $z = a/2 + ib$ is given by
%-------------------
\begin{equation}\label{eq:zhalfvelocity}
\bar{u}\vert_{z=a/2 + ib}=\lim_{N\to\infty}\dfrac{\Gamma}{2ai}
\left[ 
\sum_{n=-N, n\neq 0}^{N} \cot  \dfrac{- \pi nhi}{a} \ \ -
\sum_{n=-N}^{N}\cot \dfrac{\pi (a/2 +  (b- nh)i) }{a} 
\right] ,
\end{equation}
%-------------------
where the first sum is equal to zero and one can readily verify, upon using known trigonometric 
identities and relabeling the index $n$ by $-n$, that the second sum takes the same 
form as the expression on the right-hand side of \eqref{eq:z0vel}. 
To complete the proof, one needs to show that the infinite sum in~\eqref{eq:z0vel}
converges to a finite quantity which would be the translational
velocity $U$ of the streets. \textcolor{black}{It is known that such infinite sums
are only conditionally convergent with the result depending on the truncation
criterion, see~\cite{Oneil1989}. We use a truncation criterion which consists
of  incorporating an additional `upper' and an additional `lower' vortex street
at every step of the sum as follows} 
%-------------------
\begin{equation}\label{eq:velocity}
\begin{split}
U = \lim_{N \to \infty} U_N, \qquad U_N = \dfrac{\Gamma}{2a}
\sum_{n=-N}^{N}\tanh \dfrac{ \pi(b+ nh) }{a} .
%\\
 %\textcolor{black}{ =   \dfrac{\Gamma}{2a}\left[
%\sum_{n=-N+1}^{N-1}\tanh \dfrac{ \pi(b+ nh) }{a} + \tanh \dfrac{ \pi(b+ Nh) }{a} + \tanh \dfrac{ \pi(b - Nh) }{a}.\right]}
\end{split}
\end{equation}
%------------------
\textcolor{black}{The sum in~\eqref{eq:velocity} can be rewritten as} 
%-------------------
\begin{equation}\label{eq:velocity2}
\begin{split}
 \textcolor{black}{ U =  \dfrac{\Gamma}{2a} \tanh \dfrac{\pi b}{a} + 
\dfrac{\Gamma}{2a} \sum_{n=1}^{\infty} \left[\tanh \dfrac{ \pi(b+ nh) }{a} + \tanh \dfrac{ \pi(b - nh) }{a} \right]} .
\end{split}
\end{equation}
%------------------
\textcolor{black}{%which implies, given that the $\tanh$ function is odd, that 
%the terms generated between two consecutive steps cancel in the sum and hence 
One can readily verify -- using a standard convergence test such as the integral test (see~\cite[Section 3.3]{Knopp1956}) 
and the fact that `$\tanh$' is an odd monotonic function -- 
that the sum in~\eqref{eq:velocity2} converges. 
The rate of convergence of $U_N$ as a function of $(2N+1)$
is shown numerically in Figure~\ref{fig:velocity} 
for $b = 0.2805$ and different values of $h$.
Note that the summation criterion in~\eqref{eq:velocity} is analogous
to that employed in~\cite{Saffman1992} for the case of a single row of vortices
but differs from the one in~\cite{Oneil1989,StAr1999} for the doubly periodic lattice
which is based on radially symmetric truncations that guarantee
the terms generated between two consecutive steps in the sum to cancel out. 
The two sums, the one employed here %in~\eqref{eq:velocity}
and the one used in~\cite{Oneil1989,StAr1999}, converge to
the same value modulo a constant equal to $\Gamma b/ah$.
The values of the equilibrium velocity obtained using the Weierstrass zeta-function 
of~\cite{Oneil1989,StAr1999} plus the constant term $\Gamma b/ah$ 
are plotted on Figure~\ref{fig:velocity} for comparison purposes.
The constant term can be traced to the fact that the velocity 
of the periodic lattice in~\cite{Oneil1989,StAr1999}
is invariant to periodic shifts in both the $x$ and $y$-directions
(see~\cite[pages 103--105]{StAr1999})
while the velocity in~\eqref{eq:velocity} is invariant  to periodic shifts only in the $x$-direction.
In other words, the expressions for the stream function in~\eqref{eq:infinitepotential} and the vortex 
velocity in~\eqref{eq:velocity}, while consistent with each other, are not invariant to a coordinate 
transformation that interchanges the $x$ and $y$-directions. 
%In order to obtain such invariance, one would
%have to modify the stream function~\eqref{eq:infinitepotential} (by adding an infinite sum in the $y$-direction) 
%that guarantees that the effect of two vortices infinitely far apart in the $y$-direction on each other is zero.
}

Figure~\ref{fig:far}  shows the streamlines for multiple vortex streets. 
The parameter values are set to  $a=1$, $b=\frac{a}{\pi}\cosh^{-1}\left(\sqrt{2}\right)$, $h=1.2$, 
and $\Gamma=1$. 
The potential function in~\eqref{eq:infinitepotential} and the translational velocity in~\eqref{eq:velocity} 
are computed using $N =150$, that is, using a total of $301$ streets.
Figure~\ref{fig:far}(a) 
depicts the streamlines in inertial frame given by the contour plots of Im$(w_{N = 150})$, the imaginary part
of $w_{N = 150}$, while Figure~\ref{fig:far}(b) 
shows the streamlines in a frame moving with the translational velocity $U_{N=150}$ of
the streets which are obtained from Im$(w_{N = 150} - z  U_{N=150})$. Note that the streamline
topology in Figure~\ref{fig:far}(b) is similar to that in Figure~\ref{fig:singlestreet}(b) which means
that, for this value of $h$, the streets are not in communication with each other and each street behaves qualitatively (but not necessarily quantitatively) 
as if it were a single street. We also note that since we are using a finite sum approximation $U_N$, the 
streamline patterns begin to `distort' near the outer edges (not shown).

%%%%%%%%%%%%%%%%%%%%%%%%%%%%%%%%%%%%%%%%%%%%%%%%%%%%%

%-------------------
 \begin{figure}
 \centering
   \includegraphics[scale = 0.5] {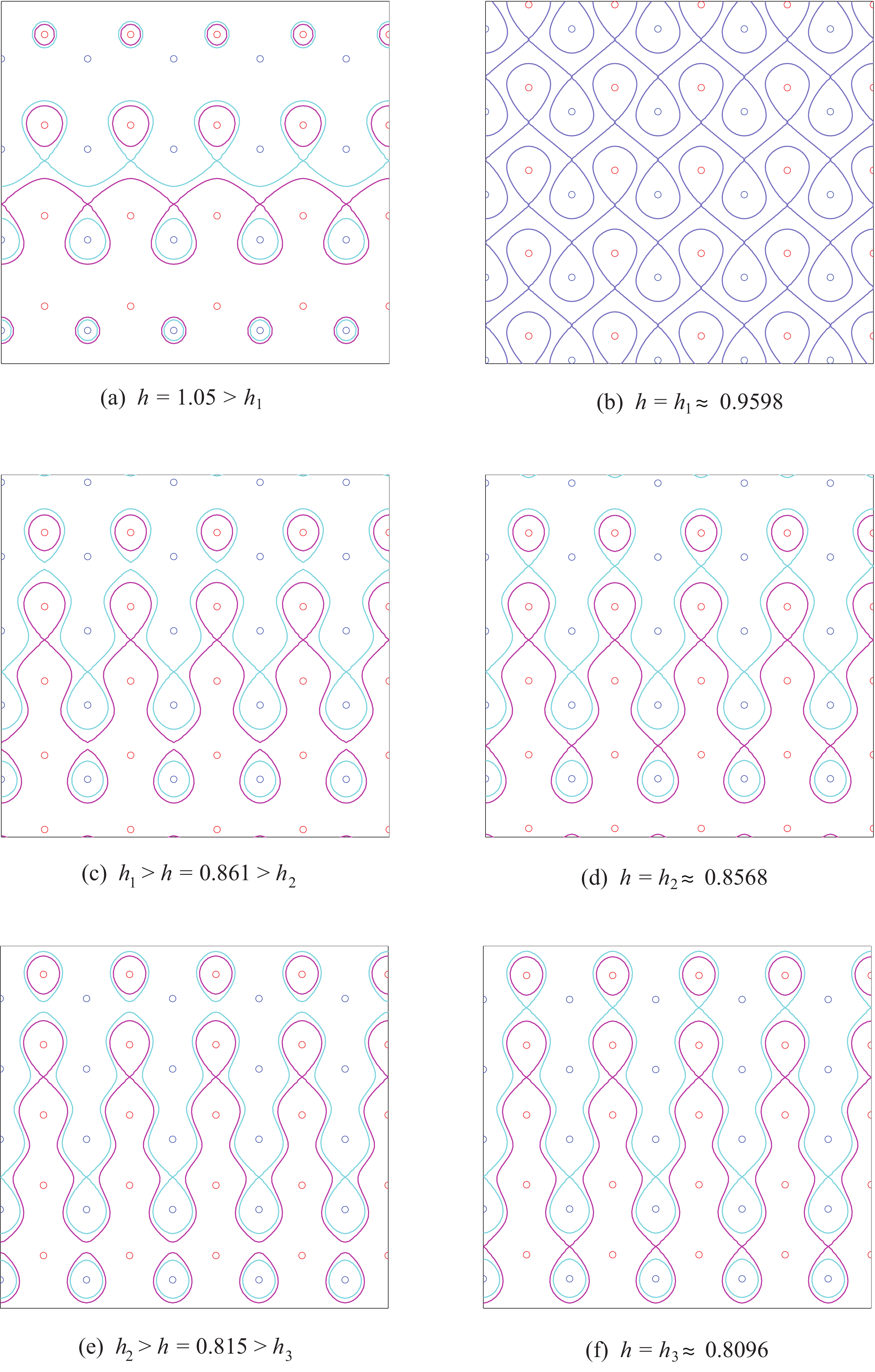}

   \caption{\textit{Streamlines plotted in moving frame for different values of $h$.
 The parameters $a$ and $\Gamma$ are normalized to unity. The parameter $b$ is set to $b=\frac{1}{\pi}\cosh^{-1}\left(\sqrt{2}\right)\approx0.2805$. The infinite sums in~\eqref{eq:infinitepotential} and~\eqref{eq:velocity} 
 are truncated using $N=150$, i.e., using a total of $301$ streets.}}
\label{fig:topology}
\end{figure}
%-------------------

%%%%%%%%%%%%%%%%%%%%%%%%%%%%%%%%%%%%%%%%%%%%%%%%%%

\section{Streamline topology}
\label{sec:stream}

The goal of this section is to study the effect of the separation distance $h$ on the streamline topology.
In particular, we identify streamline topologies that support global fluid transport between adjacent streets (as shown in
Figure~\ref{fig:far}) and others that link together two or more adjacent streets, hence supporting fluid transport
among the linked streets. 

It is convenient for this parametric study 
to rescale space and time so that $a$ and $\Gamma$ are normalized to unity.
The behavior of the system depends then on  two non-dimensional parameters, 
namely, $b$ and $h$
(which correspond to $b/a$ and $h/a$ in the dimensional system).
In Figure~\ref{fig:topology}, we fix $b$ at $b= \frac{1}{\pi}\cosh^{-1}\left(\sqrt{2}\right)$
and vary $h$. We examine the streamline topology in a frame moving
with the same translational velocity as the streets themselves.
The separatrices, that is, the streamlines associated with stagnation points
of hyperbolic type, play a central role in defining regions with distinct fluid behavior.
By stagnation points we mean points that are moving with the same translational 
velocity as the vortices themselves. Such points are found numerically
by solving for values of $z$ for which the velocity in~\eqref{eq:complexvelocity} is equal
to the translational velocity of the vortices obtained from~\eqref{eq:velocity} (for a truncation value $N = 150$). 

%-------------------
 \begin{figure}
   \centering 
   \includegraphics[scale = 0.5] {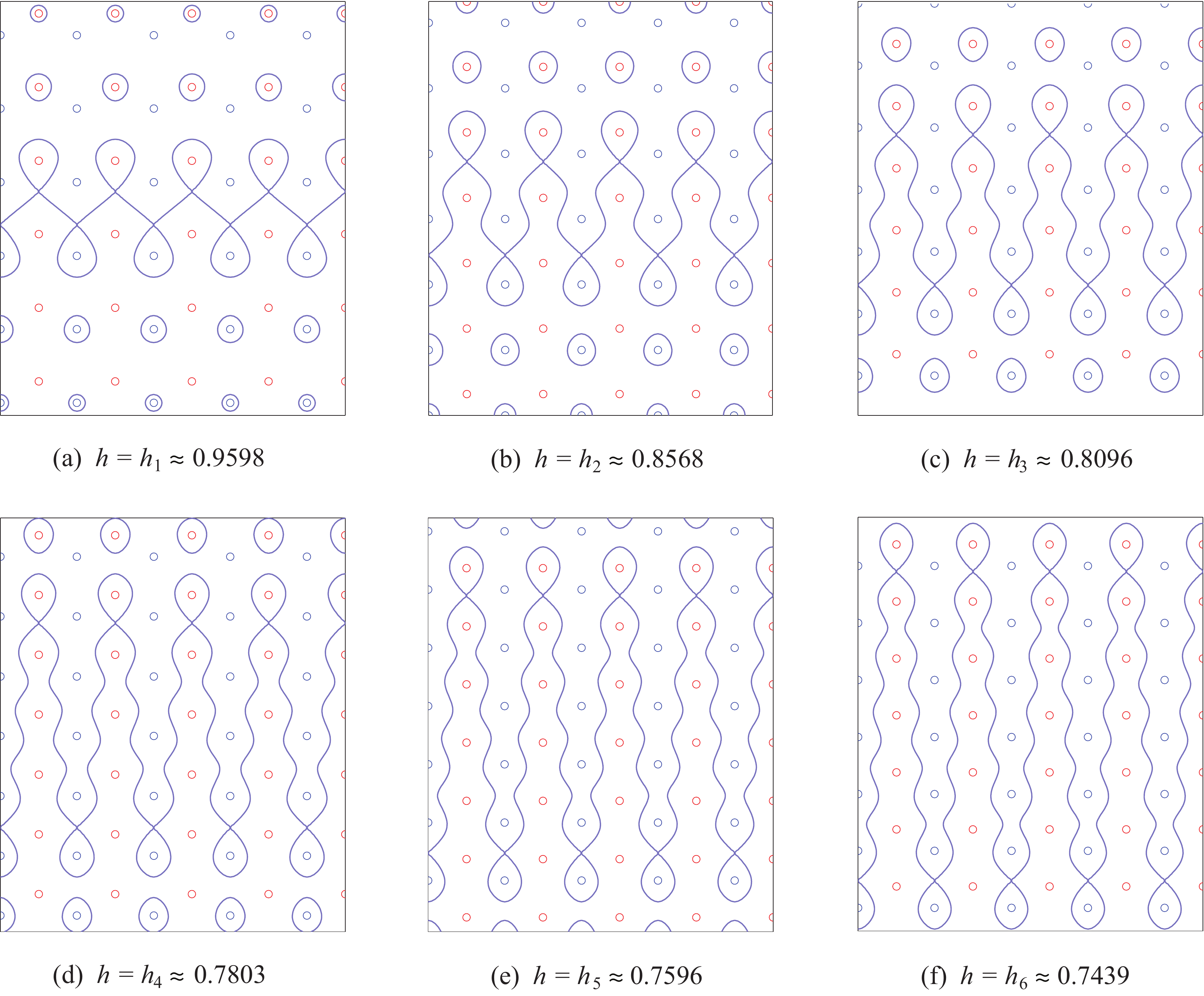}

   \caption{\textit{The first six bifurcations in the bifurcation sequence obtained by decreasing $h$.
The parameter $b$ is set to $b=\frac{1}{\pi}\cosh^{-1}\left(\sqrt{2}\right)\approx0.2805$. The infinite sums in~\eqref{eq:infinitepotential} and~\eqref{eq:velocity} 
 are truncated using $N=150$, i.e., using a total of $301$ streets.}}
\label{fig:bifurcation}
\end{figure}
%-------------------  

The streamlines shown in Figure~\ref{fig:topology}(a) correspond to $h=1.05$.
Only a small number of streamlines are shown in order to emphasize those that 
pass through the stagnation points, thus forming separatrices 
between regions of distinct fluid behavior. 
For this value of $h$, one can identify three regions
of distinct fluid behavior: (1) around each point vortex, there is a region containing the 
vortex and bounded by a separatix where the fluid orbits around the vortex; (2) within each street,
there is a region bounded by two separatrices and void of point vortices where the 
fluid is transported globally while touring around opposite-sign vortices of the same street; and (3)
between two adjacent streets, there is a region where the fluid is transported globally.
One can readily observe that  the streamlines in Figure~\ref{fig:topology}(a)   
have the same topology as those depicted in Figure~\ref{fig:far}(b) for $h=1.2$ but, 
in Figure~\ref{fig:topology}(a), the separatices have moved closer to each other making 
the region of global transport between two adjacent streets smaller. 
Upon decreasing $h$, the streamline topology remains the same  until $h$ 
reaches the critical value $h_1 \approx0.9598$, shown in Figure~\ref{fig:topology}(b). 
At this bifurcation value, the separatrix associated with a given street  collapses onto 
the separatrix associated with its immediate neighbor (or adjacent street). 
As a result, the region of global fluid transport between two adjacent streets disappears. 
For $h< h_1$, one has a different streamline topology  consisting of two regions 
of distinct fluid behavior as shown in Figure~\ref{fig:topology}(c). Namely,
there is a region around each vortex containing the 
vortex and bounded by a separatix where the fluid orbits around that vortex
and there is a second region linking two adjacent streets where the fluid
is transported globally along paths that alternate around opposite-sign vortices of two 
adjacent streets. This streamline topology remains the same until $h$ reaches a second 
critical value $h_2 \approx0.8568$, shown in Figure~\ref{fig:topology}(d).

At this bifurcation value the separatrix linking two adjacent streets collapses 
onto the neighboring separatrix (also linking two adjacent streets) thus creating 
one separatix linking three adjacent streets. For $h< h_2$, one
has a different streamline topology  consisting of two regions of distinct
fluid behavior, see Figure~\ref{fig:topology}(e). Namely, there is
a region around each vortex where the fluid orbits around the vortex 
and a second region linking three adjacent streets where the fluid
is transported globally along paths alternating around opposite-sign vortices of the 
two non-adjacent linked streets. This streamline topology is maintained until
a third bifurcation occurs at $h_3 \approx 0.8096$ where 
two adjacent separatices merge again hence linking four adjacent
streets. This pattern  gets repeated when decreasing $h$ further, thus 
creating a sequence of bifurcations where, at each bifurcation, 
one additional street joins the streets that are already linked.
Figure~\ref{fig:bifurcation} shows the separatrices associated
with the first six bifurcations.

%%%%%%%%%%%%%%%%%%%%%%%%%%%%%%%%%%%%%%%%%%%%%%%%%%%%%
\section{Bifurcation analysis}
\label{sec:bif}

\begin{figure}
   \centering
   \includegraphics[scale=0.5,angle=0] {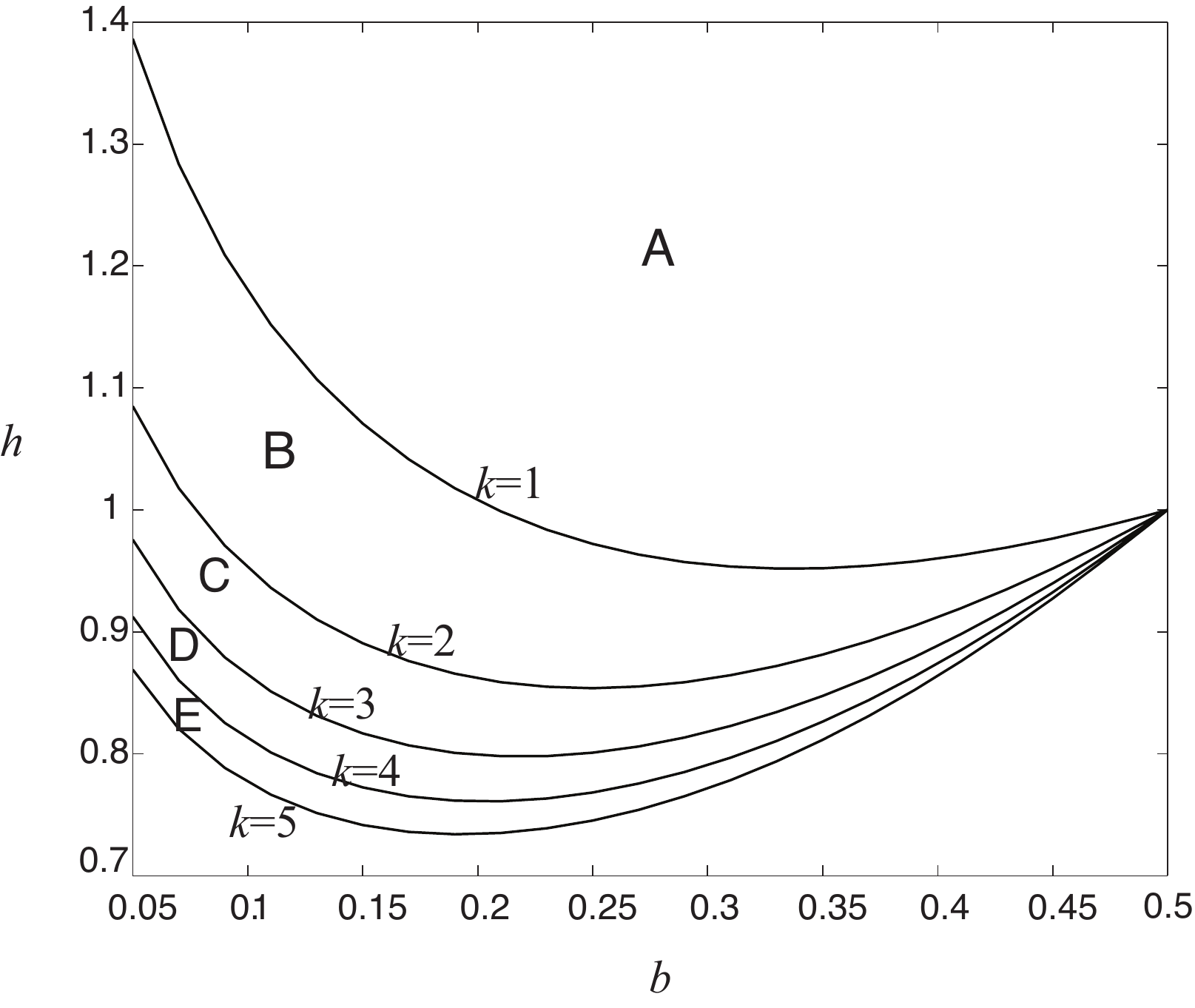}

   \caption{\textit{Bifurcation curves in the parameter space $h$ and $b$. The parameter
   $b$ is varied from $0.05$ to $0.5$ by increments of $0.45/22$. For each value of $b$
   the first five bifurcation values of $h$ are computed. The obtained curves separate
   the parameter space into regions that correspond to distinct streamline topologies. 
   For example, region $B$ is characterized by streamlines that support communication 
   and fluid exchange between two neighboring streets while, in region $C$, 
   fluid is transported between three neighboring streets, and so on. }}
\label{fig:BifCurves}
\end{figure}

\begin{figure}
   \centering
   \includegraphics[scale=0.5,angle=0] {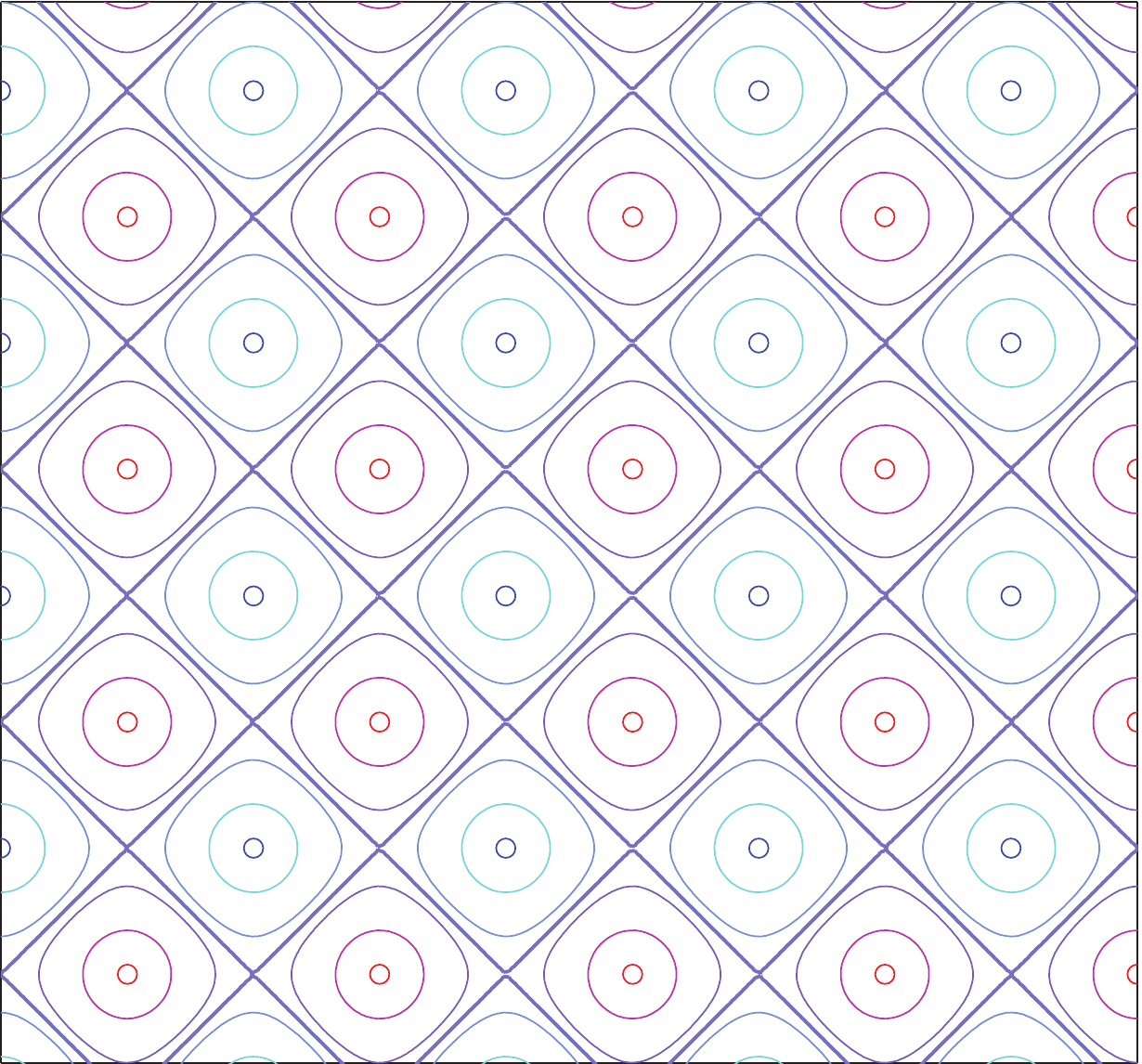}

   \caption{\textit{Streamlines for $b=0.5$, $h=1$ where all bifurcation curves in Figure~\ref{fig:BifCurves} intersect.
   The resulting lattice is an `Abrikosov'   square lattice (see Abrikosov (2004)).}}
\label{fig:BifIntersect}
\end{figure}

The bifurcation sequence shown in Figure~\ref{fig:bifurcation} is persistent when varying 
the parameter $b$ as shown by the bifurcation curves in Figure~\ref{fig:BifCurves}. 
These curves are obtained by varying $b$ from $0.05$ to $0.5$ using discrete increments 
and, for each value of $b$,  computing the bifurcation values of $h$ as done 
in~Section~\ref{sec:stream}. Only the curves associated with the first five bifurcations are shown,
which define five regions in the parameter space labeled A to E.
For parameter values $(h,b) \in A$, the streets are decoupled and fluid is transported globally
between streets (similarly to the streamline topology shown in Figure~\ref{fig:topology}(a)) 
while for $(h,b) \in B$,  two adjacent streets get coupled in the sense 
that the fluid transported globally \textcolor{black}{moves} along paths that alternate between two adjacent streets
(similarly to the streamline \textcolor{black}{topology} shown in Figure~\ref{fig:topology}(c)), and so on.
Note that all bifurcation curves intersect at the point $b = 0.5$ and $h=1$, which
corresponds to a  square `Abrikosov'  lattice (Abrikosov (2004)) as shown in Figure~\ref{fig:BifIntersect}.
That is, the bifurcation sequence reduces to one bifurcation point where all the separatrices
collapse onto one. By symmetry, this square lattice is stationary.

%-------------------
%-------------------

%-------------------
%-------------------

\begin{figure}
   \centering
   \includegraphics[scale = 0.5,angle=0] {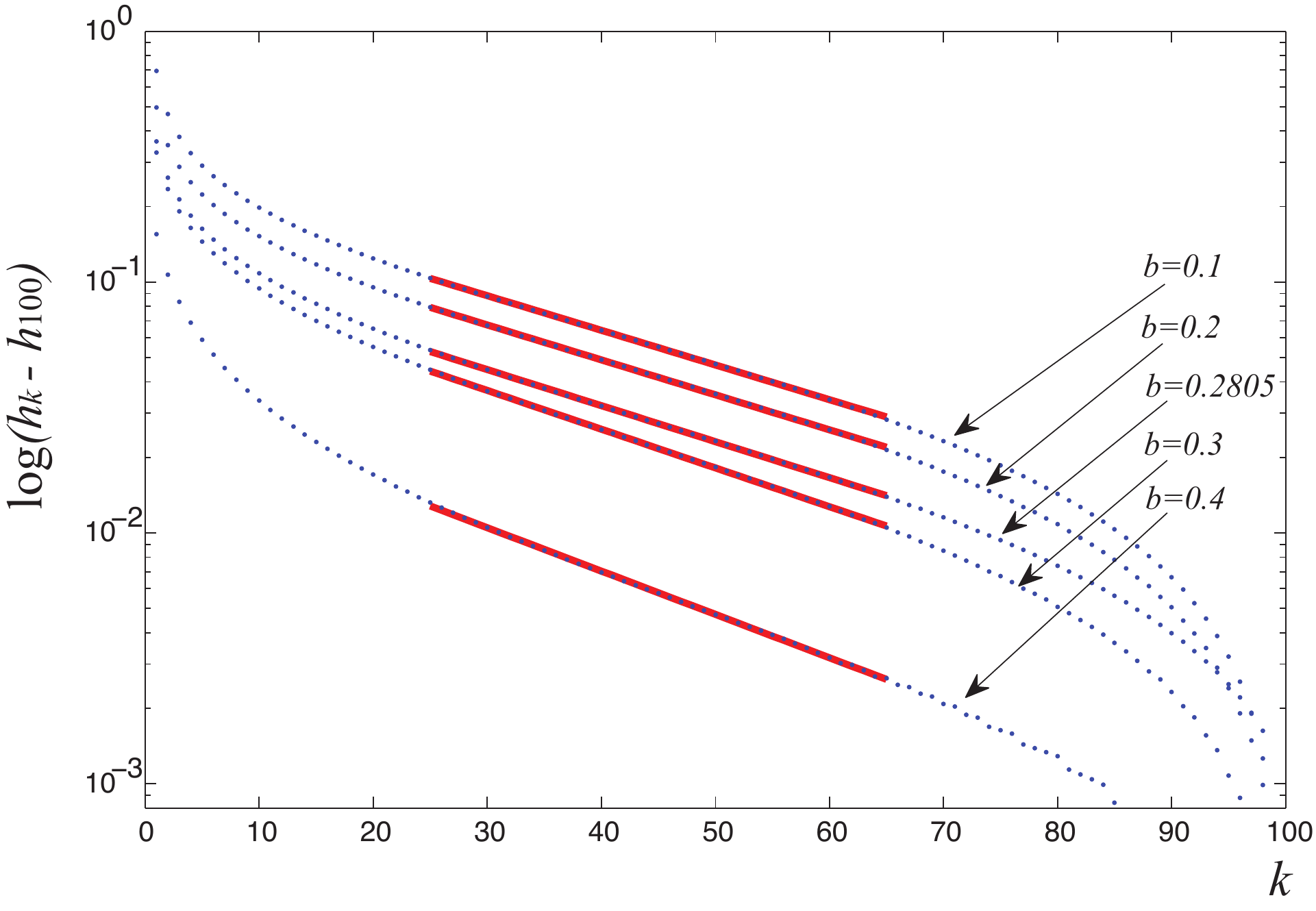}
   
   \caption{\textit{The values of $\log{\left(h_k-h_{100}\right)}$ versus the bifurcation index $k$
    for five distinct values of $b$, namely, $b=0.1, 0.2, 0.2805, 0.3, 0.4$.}}
\label{fig:fit}
\end{figure}

We conjecture that, for a given $b$ ($0<b<0.5$), there exists a countably  infinite number of bifurcation
values  $h_k$, where $k$ is the bifurcation index $k = 1,2,\ldots$, and we  postulate the following scaling law for $h_k$,
%-------------------
\begin{equation}\label{eq:scaling}
h_k \sim h_{\infty}+\dfrac{c}{\delta^k}, 
\end{equation}
%-------------------
for $k$ {\it sufficiently} large, and
where $h_{\infty}$ is the bifurcation value as $k\to \infty$ while $c$ and $\delta$ are  to be determined. 
\textcolor{black}{Upon taking the logarithm of}~\eqref{eq:scaling} \textcolor{black}{one gets}
%-------------------
\begin{equation}\label{eq:scaling2}
\log(h_k -h_{\infty}) \sim  -k\log\delta  + \log c.
\end{equation}
%-------------------
That is, the value $-\log{\delta}$ is the slope of the curve obtained by plotting 
$\log{\left(h_k-h_{\infty}\right)}$ versus $k$ in a half $\log$ plane. 
Since we do not have a closed form expression for the bifurcation $h$
and, therefore, cannot compute the value of $h_{\infty}$ analytically,
we use   $h_{100}$ as an approximation for $h_\infty$. That is, we compute
the first $100$ bifurcations and set $h_\infty = h_{100}$. As before, the
value of $N$ is set to $N=150$, that is, the infinite array of streets
is discretized using $301$ streets.
The values of $h_{100}$, $c$ and $\delta$ are determined numerically for 
$b = 0.1, 0.2, 0.2805, 0.3$ and $0.4$ and are reported in Table~\ref{tbl:bifs}. 
A half-logarithmic plot  of $\log{\left(h_k-h_{100}\right)}$ versus the bifurcation index $k$
is shown in Figure~\ref{fig:fit} for all $100$ bifurcation values of $h$.
Superimposed on these, we plot the right hand side of~\eqref{eq:scaling2}, 
namely, $-k (log\delta) + log c $, which are the red lines in Figure~\ref{fig:fit}.
Clearly, for each value of $b$, there is an interval where the scaling law
fits the computed bifurcations, which confirms that the bifurcation parameter $h$
scales as suggested in~(\ref{eq:scaling},\ref{eq:scaling2}).

%--------------------------------------------------------------
\begin{table}
\begin{center}
\begin{tabular}{||c|c|c|c||}
    \hline
    \hline
    $b$ & $h_{100}$ & $c$ & $\delta$ \\ 
    \hline 
$0.1$ & $0.4852$  & ${0.2291}$ & ${1.0323}$ \\
$0.2$ &  $0.5110$  & ${0.1766}$ & ${1.0326}$ \\
$0.2805$ & $0.5963$ & ${0.1206}$ & ${1.0336}$ \\
$0.3$ & $0.6267$  & ${0.1074}$ & ${1.0362}$ \\
$0.4$ & $ 0.8051$ & ${0.0346}$ & ${1.0406}$ \\
    \hline 
    \hline
    \end{tabular}
    \end{center}
\vspace{-.4cm}
    \caption{\textit{The scaling law is given by $h_k =h_{\infty}+\dfrac{c}{\delta^k}$ where $h_\infty$ ia approximated by $h_{100}$. The values of $h_{100}$, $c$ and $\delta$ are determined for $b = 0.1, 0.2, 0.2805, 0.3$ and $0.4$. A half-logarithmic plot 
    of the scaling law for these parameter values is shown in Figure~\ref{fig:fit}.}}
\label{tbl:bifs}
    \end{table}
%--------------------------------------------------------------
%-------------------
%-------------------

Since the scaling law is postulated for values of $k$ sufficiently large, we do not expect the
curve to be perfectly linear in the {\it small} $k$ regime (which it is not). In addition, as $k$ approaches
$100$, our use of $h_{100}$ instead of $h_\infty$ as well as the finiteness of our lattice in the $y$ direction
destroys the linear scaling.

%%%%%%%%%%%%%%%%%%%%%%%%%%%%%%%%%%%%%%%%%%%%%%%%%%%%%

\section{Time Evolution of the Streamlines}
\label{sec:timescales}

The prevalence of the streamline topologies reported in Sections~\ref{sec:stream} and~\ref{sec:bif} 
when subject to spatial and temporal perturbations is a central issue in arguing their validity. 
In this section, we do not pursue a rigorous stability
analysis of the relative equilibria for the infinite array of vortex streets. Such analysis would require 
mathematical tools similar to those employed in
\cite{Stremler2003,StAr1999} and~\cite{Tkachenko1966a,Tkachenko1966b} and 
would be restricted to a class of periodic perturbations 
of the infinite array of streets. Instead, we investigate the dynamic evolution of a finite approximation
of the infinite array where the streets are truncated both vertically and  horizontally,
reflecting a finite number of fish within the school whose wake is diminished by the 
effects of viscosity. The truncated configuration can be viewed as a non-periodic perturbation
of the infinite one and its time evolution does provide insight into the mechanism by which the
street configuration gets destroyed. For short times, however, we argue that the analysis carried 
out for the infinite lattice case does indeed set the pattern that should be relevant on some 
timescale for the finite approximation. % long enough to be relevant in the case of a finite collection of fish.

%-------------------
 \begin{figure}
   \centering 
   \subfigure[Infinite array of infinite streets]
   {\includegraphics[height=3.8cm,angle=0] {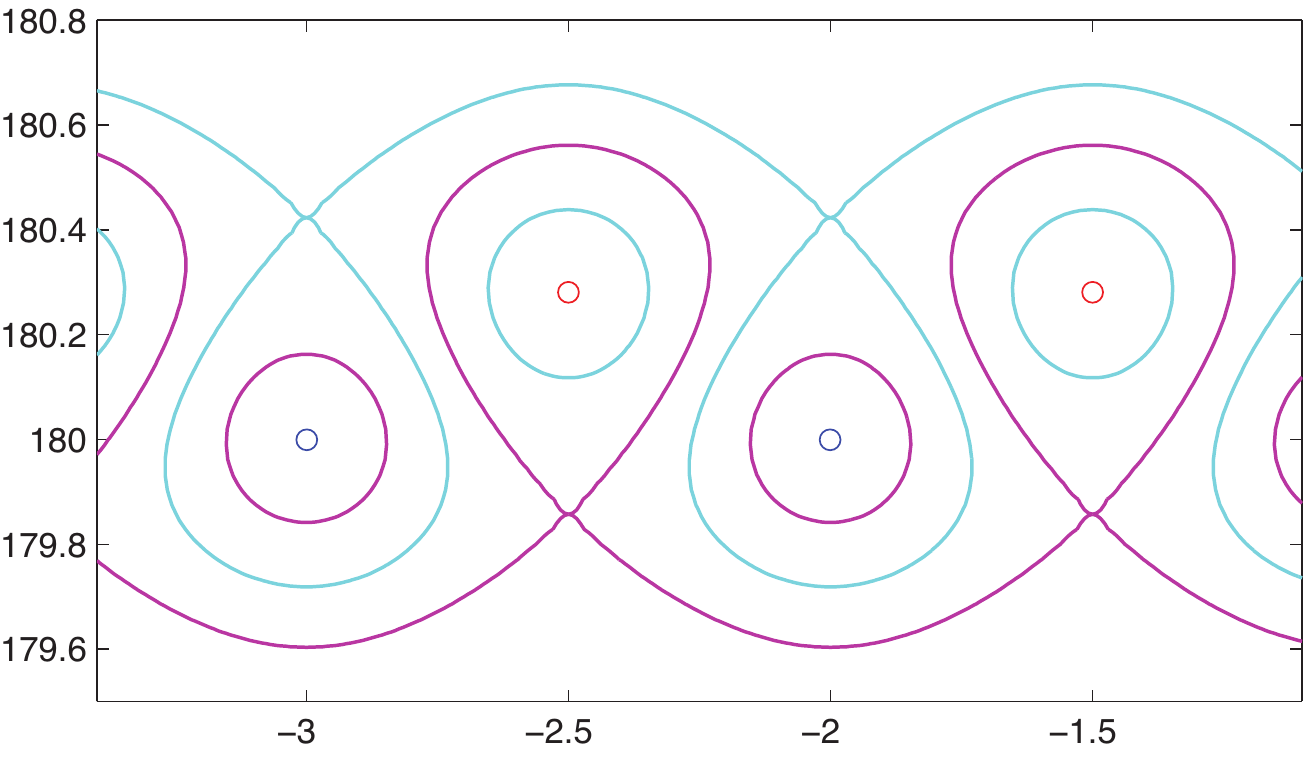}}
   \hspace{.2cm}   
   \subfigure[Finite array of truncated streets]
   {\includegraphics[height=3.8cm,angle=0] {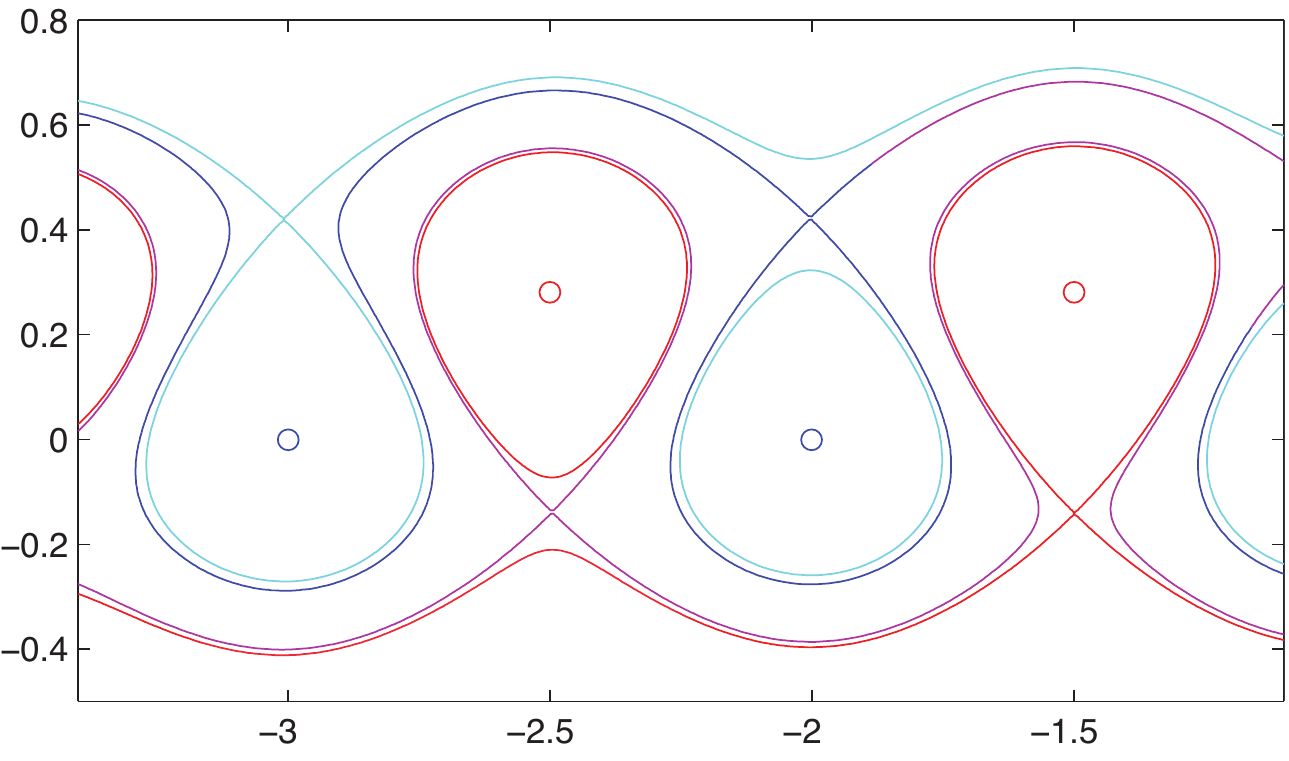}}
   
   \caption{ \textit{ Streamline topologies at time $t=0$ for 
   parameter values $a=1$, $b=0.2805$, $h=1.2$ and $\Gamma=1$. (a) Streamline plots shown in a relative 
   frame for a lattice including 301 number of infinite streets. The translational velocity of 
   the relative frame is given by~\eqref{eq:velocity}. Only a small window is shown depicting a portion of the 
   middle street out of a total of 301 streets. (b) Streamline plots shown in a relative frame for a finite array 
   consisting of a total of 11 streets each containing 2 rows of 10 point vortices (20 point vortices in each street). 
   The streamlines are shown in a moving frame whose translational velocity is taken to be equal to the 
   instantaneous velocity of the point vortex at the origin. Unlike in the infinite case where all vortices 
   have the same translational velocity, the point vortices in the shown window in (b) have velocities 
   that are slightly different from the velocity of the vortex at the origin. Most importantly, in the finite case, 
   the streamlines \textcolor{black}{corresponding to the hyperbolic stagnation} points located within the street do not coincide. 
   In other words, the hetro-clinic streamlines (or separatrices) observed in (a) have turned into a number 
   of homo-clinic separatrices in (b). }}
   \label{fig:finitevsinfinite}
\end{figure}
%-------------------

Consider a finite number of $2N+1$ truncated vortex streets where each street includes $2$ rows of finitely-many vortices  of equal and opposite strength initially placed in the street configuration described in 
Section~\ref{sec:formulate}. It is convenient for the purposes of this section to label the point vortices 
by $z_{pn}$ where, at time $t = 0$, $z_{pn} = m a + nhi$ for $p= 2m+1$ odd and $z_{pn} = (m+1/2)a + i(b+nh)$ 
for $p = 2m$ even. Also, the strength $\Gamma_{pn}$ of a point vortex labeled by $z_{pn}$ is 
given by $(-1)^p\Gamma$. Clearly, the finite configuration does not form a relative equilibrium. 
That is, starting with this `street' configuration, the vortices interact dynamically in time and their 
evolution is such that the street configuration gets broken.  Before analyzing this time evolution, 
it is informative to examine its streamline topology at time $t=0$ and 
compare that to the case of the infinite arrays considered earlier. In Figure~\ref{fig:finitevsinfinite} 
we show the difference between the streamline topologies at time $t=0$ of the infinite and finite 
arrays for parameter values $a=1$, $b=0.2805$, $h=1.2$ and $\Gamma=1$.  In Figure~\ref{fig:finitevsinfinite}(a), 
the streamlines of an infinite array of streets are plotted in a frame moving with a translational velocity 
given by~\eqref{eq:velocity} using an approximation of $301$ vortex streets, exactly as done in Figure~\ref{fig:far}. 
In Figure~\ref{fig:finitevsinfinite}(a), only a small window is shown depicting a portion of the middle street. 
Figure~\ref{fig:finitevsinfinite}(b) shows the streamlines associated with a finite array consisting of only 
11 streets each containing 2 rows of 10 point vortices each (i.e., a total of 20 point vortices in each truncated street). 
The streamlines are shown in a moving frame whose translational velocity is taken to be equal to the instantaneous velocity of the point vortex at the origin.  The (complex conjugate of the) velocity of a point vortex identified by its position $z_{kl}$ is given by
%---------
\begin{equation}\label{eq:ev}
\color{black}  \dfrac{d \bar{z}_{kl}}{dt} = \sum_{p,q,~z_{pq}\neq z_{kl}}  -   \dfrac{i\Gamma}{2\pi} \dfrac{(-1)^{p}}{z_{pq}-z_{kl}}
\end{equation}
%--------
Unlike in the infinite case where all vortices have the same translational velocity,
the point vortices in the finite case have different velocities. Also in the finite case, 
the streamlines corresponding to the hyperbolic stagnation points located within the street do not coincide. 
In other words, the hetro-clinic streamlines or separatrices observed in Figures~\ref{fig:far}\textcolor{black}{(b)} and~\ref{fig:finitevsinfinite}(a) for the infinite case turn into a number of closely-spaced homo-clinic separatrices in  Figure~\ref{fig:finitevsinfinite}(b) creating new 
regions of fluid transport between these separatrices not present in the infinite case while maintaining 
the three regions of fluid transport identified in Figure~\ref{fig:far} for the infinite array. 
\textcolor{black}{This splitting of the hetero-clinic streamlines is analogous
to the splitting of separatrices under small perturbations, see for example~\cite{GeLa2001}.}
We repeated this exercise for a range of values of the parameter $h$ corresponding 
to various streamlines topologies 
and observed similar results (not shown). 
\textcolor{black}{Namely, at time $t=0$, instead of the hetero-clinic streamlines  identified in Section~\ref{sec:stream} and defining the separation between the various regions of fluid transport, 
in the finite approximation, 
one obtains a family of closely-spaced homo-clinic streamlines defining the separation 
between the same regions.
In this sense, we say that the general streamline patterns in the infinite array prevail
in the finite approximation.}

\begin{figure}
\centering 
\includegraphics[scale = 0.5] {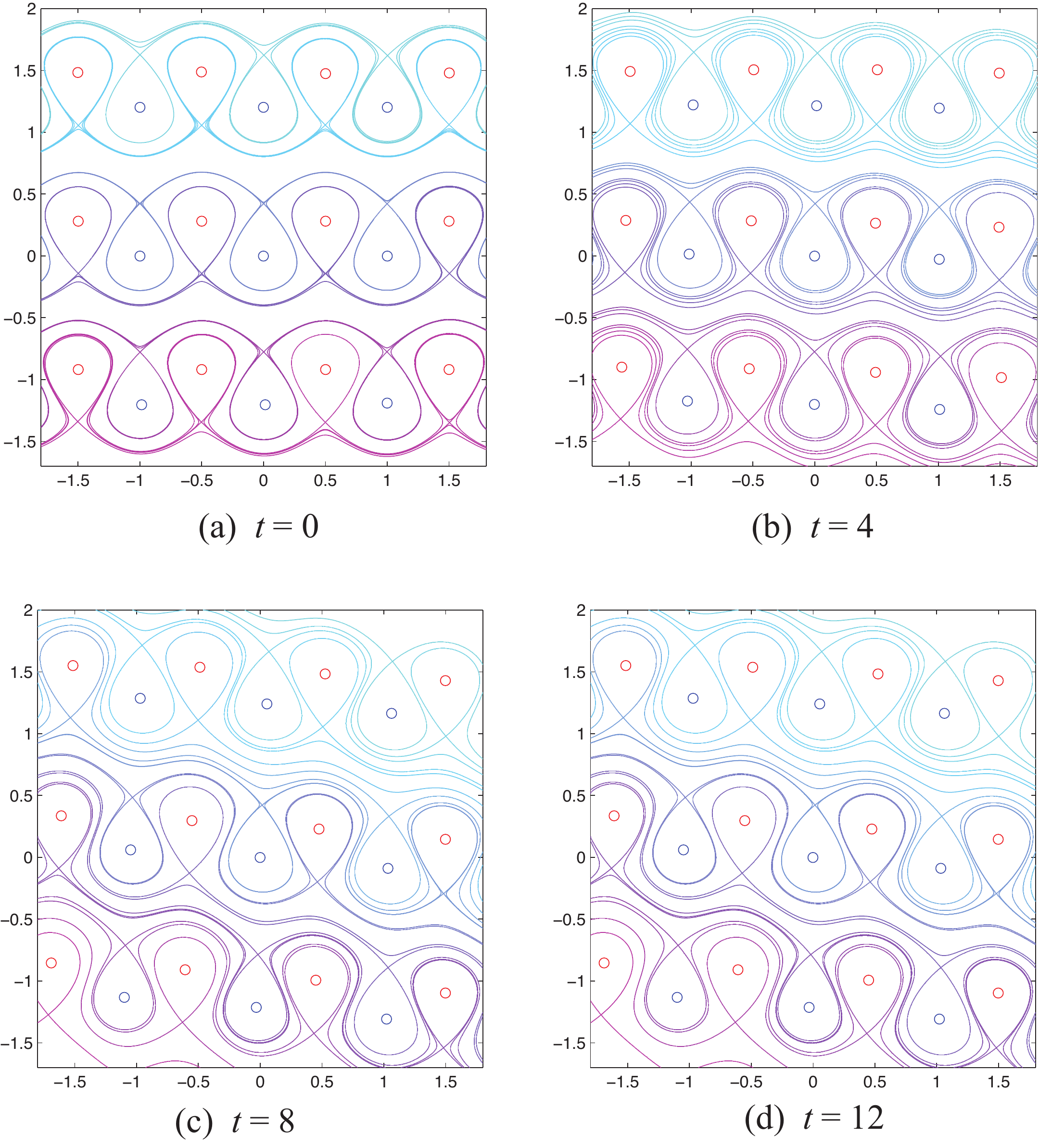}

\caption{\textit{Snapshots of the streamlines plotted in a frame moving with the same velocity as the point vortex  initially placed in the middle of the finite array of streets. This array includes 11 streets, each composed of 20 equal and opposite point vortices  (yielding a total of 220 point vortices). The parameters that define the initial spacing of the vortices are set to $a=1$, $b=0.2805$, $h=1.2$ while the vortex strength is normalized to 
 $\Gamma=1$. Clearly, the initial street-like configuration changes with time and eventually breaks completely.}}
\label{fig:evolutionTopologyA}
\end{figure}

\begin{figure}
\centering 
\includegraphics[scale = 0.5] {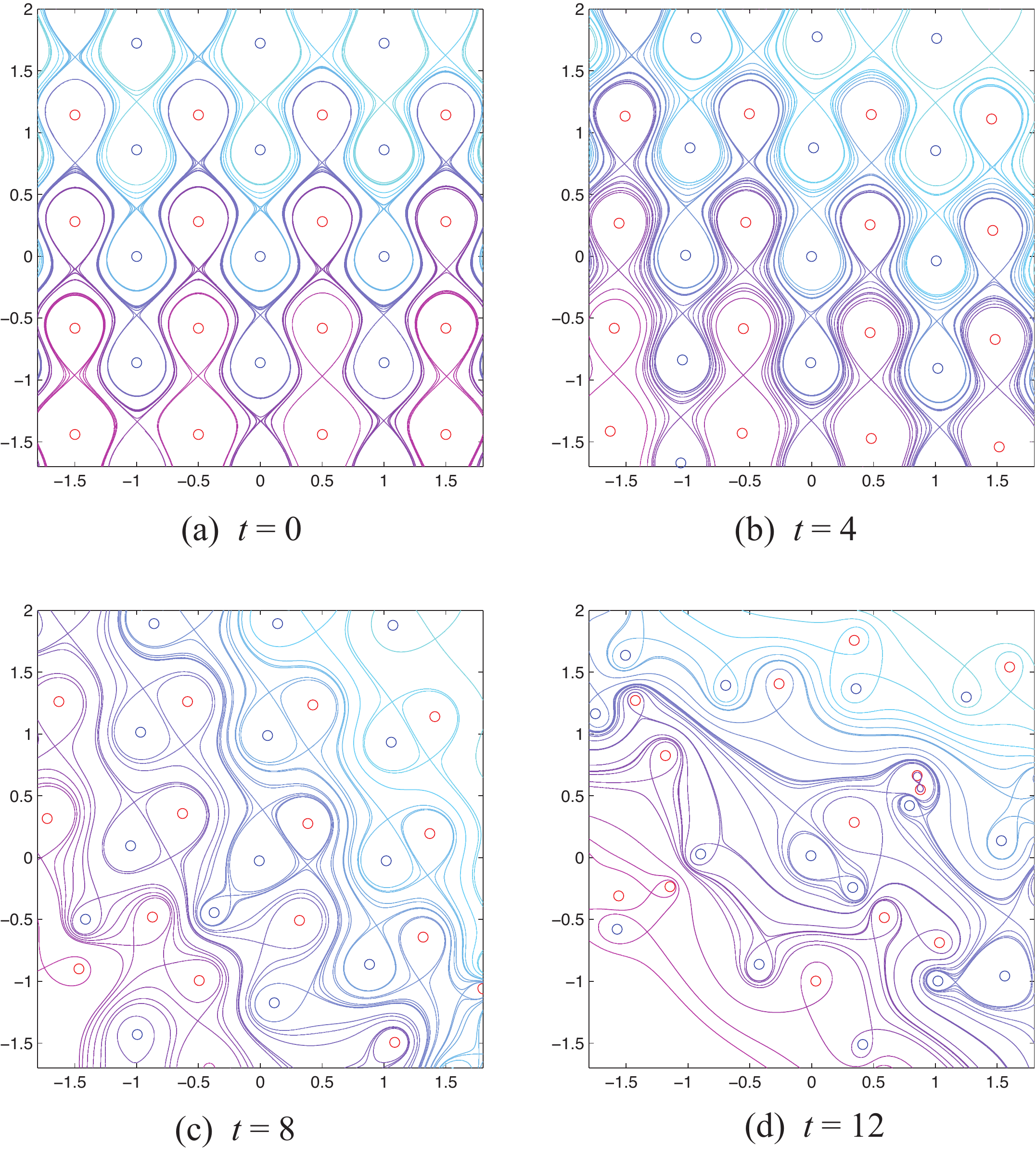}

\caption{\textit{Snapshots of the streamlines plotted in a frame moving with the same velocity as the point vortex  initially placed in the middle of the finite array of streets. This array includes 11 streets, each composed of 20 equal and opposite point vortices  (yielding a total of 220 point vortices). The parameters that define the initial spacing of the vortices are set to $a=1$, $b=0.2805$, $h=0.861$ while the vortex strength is normalized to 
 $\Gamma=1$. Clearly, the initial street-like configuration changes with time and eventually breaks completely.}}
\label{fig:evolutionTopologyB}
\end{figure}

We use~\eqref{eq:ev} to describe the time evolution of a finite array consisting of 11 streets, each composed of $20$ point vortices of equal and opposite circulation, thus leading to a total of $220$ point vortices and $440$ coupled nonlinear ODEs governing the positions of the point vortices, namely $x_{kl}$ and $y_{kl}$. This set of equations can be numerically integrated to show the evolution of the vortex array.
In Figures \ref{fig:evolutionTopologyA} and \ref{fig:evolutionTopologyB} we show typical evolution diagrams
of the finite streets at three successive snapshots ($t = 4, 8, 12$) for two different values of $h$, all other parameters equal. As expected, because the street is finite and thus not an exact equilibrium configuration, 
the clean streamline
patterns and scaling shown for the infinite street cases \textcolor{black}{break} down at some finite time.
As shown in Figure \ref{fig:evolutionTopologyA} ($h = 1.2$), the topological pattern still holds at $t = 4$, 
but begins to break down in an interesting way by the time $t = 8$. By the time $t = 12$, the system clearly 
breaks into clusters of
translating dipole pairs of equal and opposite groups of point vortices translating at an oblique angle. 
In Figure \ref{fig:evolutionTopologyB} ($h = 0.861$), the prevelant global topological patterns discussed in the inifinite
street case still seems to hold at $t = 4$, but has broken down by the time $t = 12$. Here, since the streets are
closer, the breakdown occurs more quickly, showing that by the time $t = 12$, the clusters no longer remain as 
dipole pairs, but group into more complex combinations of co-rotating pairs and tripoles. 
\textcolor{black}{The streamline patterns from the
evolution of truncated wakes are reminiscent of those seen in careful experiments of interacting 2D 
wakes, such as those reported in~\cite{RaAlPaHo1988}}.
Our two main conclusions from these finite-array simulations are that (i) the streamline topology of the finite array, for short times, can be expected to roughly follow our analysis of the previous sections, but (ii) eventually (depending on the various parameters), the topology breaks down and the finite collection of vortices disperses presumably chaotically, but certainly not in a clean von K\'{a}rm\'{a}n arrangement.

%%%%%%%%%%%%%%%%%%%%%%%%%%%%%%%%%%%%%%%%%%%%%%%%%%%%%%%%%%

%-------------------
%-------------------

%-------------------
%-------------------

\section{Conclusions}
\label{sec:conc}

The  multiple wake system developed in this paper is proposed as a model of 
the mid-wake region generated by  a school of fish, swimming in a coordinated fashion.
Immediately  behind each individual fish (i.e. the near-wake), a classical (\textcolor{black}{reverse})  von K\'{a}rm\'{a}n vortex street
would be a standard model of the flow structure, whereas far downstream (i.e. the far-wake),
 viscous effects effectively diminish the presence of the entire wake from the school. However in the mid-wake region, the wake generated from each individual fish interacts nonlinearly with its neighbors  to produce a complex structure whose streamline pattern depends on (i) the aspect ratio of each individual  von K\'{a}rm\'{a}n vortex street, and (ii) the distance between neighboring
streets. Our model shows that the closer the wakes are to their neighbors, the `wider' the cooperative effect 
becomes, in the sense that more and more neighboring wakes are able to communicate and exchange fluid. 
The structures potentially affect the  behavior of a fish swimming in this mid-wake region in two important 
ways. First, the separatrices produced by the field of wakes \textcolor{black}{may effectively block passive locomotion of fish. There is however no guarantee that these separatrices which constitute barriers to passive tracers would also constitute barriers to passive locomotion of fish who, unlike passive tracers, may alter the streamline patterns in the flow (see for example~\cite{KaOs2008}).}
%passage of a fish making it
%more difficult for the individual to re-arrange its position with respect to the group. 
Second, because of this
`blocking' effect, passive particles advected in the wake, such as nutrients, can be more effectively trapped, 
leading to benefits for the individual fish. The bifurcations of the streamline patterns were studied as a function of
the wake parameters and shown to follow a geometric scaling law, not unlike the scaling law produced in a
Fiegenbaum period-doubling cascade. Evolution patterns for finite-arrays show that the scaling laws and
global streamline patterns discussed in the infinite array setting eventually break down and become more 
complex. However, for short times \textcolor{black}{dictated by the aspect 
ratio of each street and the separation distance between the streets}, near the center of the mid-wake 
regions, because the effect of the perturbations due to the finite-array has not yet been felt, the pattern 
from the infinite array prevails and is expected to play some determining role in the collective dynamics 
of the school. Future work will focus on whether these `cooperative' features identified in 
this idealized setting persist in more realistic models which incorporate still more 
dynamical complexity \textcolor{black}{such as the presence of moving boundaries that model the fish}
and in three-dimensions.
%
%%%%%%%%%%%%%%%%%%%%%%%%%%%%%%%%%%%%%%%%%%%%%%%%%%%%%%

%\appendix
%\section{Convergence test}
%\setcounter{equation}{0}
%\def\theequation{A.\arabic{equation}}

%
%The sum in~\eqref{eq:velocity} can be rewritten as 
%%-------------------
%\begin{equation}\label{eq:velocity2}
%\begin{split}
% \textcolor{black}{ U_N =  \dfrac{\Gamma}{2a}\left[
%\sum_{n=-N+1}^{N-1}\tanh \dfrac{ \pi(b+ nh) }{a} + \tanh \dfrac{ \pi(b+ Nh) }{a} + \tanh \dfrac{ \pi(b - Nh) }{a} \right]} .
%\end{split}
%\end{equation}
%%------------------

%
%*********************************BIBLIOGRAPHY*******************************
%

\end{document}